\documentclass[11pt]{article}
\usepackage{amsmath,amsthm,amsfonts,amssymb,enumerate}
\usepackage[colorlinks=true,citecolor=black,linkcolor=black,urlcolor=blue]{hyperref}
\usepackage{algorithm}
\usepackage{algorithmic}   
 \floatname{algorithm}{Algorithm}

\usepackage[lmargin=30mm,rmargin=30mm,bmargin=30mm,tmargin=30mm]{geometry}
\usepackage[numbers,sort&compress]{natbib}

\renewcommand{\baselinestretch}{1.1}
\renewcommand{\thefootnote}{\fnsymbol{footnote}}	
\allowdisplaybreaks

\theoremstyle{plain}
\newtheorem{theorem}{Theorem}

\begin{document}

\title{\bf A linear-time algorithm for finding a complete graph minor
  in a dense graph}

\author{Vida Dujmovi{\'c}\,\footnotemark[4] \qquad Daniel J.\ Harvey
  \footnotemark[3]\qquad Gwena\"el Joret\,\footnotemark[2] \qquad
  \\ Bruce Reed \footnotemark[5] \qquad David~R.~Wood\,\footnotemark[1]}

\date{\today}

\maketitle

\begin{abstract}
Let $g(t)$ be the minimum number such that every graph $G$ with average degree $d(G) \geq g(t)$ contains a $K_{t}$-minor. Such a function is known to exist, as originally shown by Mader. Kostochka and Thomason independently proved that $g(t) \in \Theta(t\sqrt{\log t})$. This article shows that for all fixed $\epsilon > 0$ and fixed sufficiently large $t \geq t(\epsilon)$, if $d(G) \geq (2+\epsilon)g(t)$ then we can find this $K_{t}$-minor in linear time. This improves a previous result by Reed and Wood who gave a linear-time algorithm when $d(G) \geq 2^{t-2}$.
\end{abstract}

\footnotetext[4]{School of Mathematics and Statistics \&
Department of Systems and Computer Engineering, Carleton University,
  Ottawa, Canada (\texttt{vida@scs.carleton.ca}). Supported by NSERC
  and an Endeavour Fellowship from the Australian Government.}

\footnotetext[2]{D\'epartement d'Informatique, Universit\'e Libre de
  Bruxelles, Brussels, Belgium
  (\texttt{gjoret@ulb.ac.be}). Postdoctoral Researcher of the Fonds
  National de la Recherche Scientifique (F.R.S.--FNRS). Supported by
  an Endeavour Fellowship from the Australian Government.}

\footnotetext[3]{Department of Mathematics and Statistics, The
  University of Melbourne, Melbourne, Australia
  (\texttt{d.harvey@pgrad.unimelb.edu.au}). Supported by an Australian Postgraduate Award.}

\footnotetext[5]{Canada Research Chair in Graph Theory, supported by
  NSERC. School of Computer Science, McGill University,
  Montr\'eal, Canada (\texttt{breed@cs.mcgill.ca}); Laboratoire I3S,
  Centre National de la Recherche Scientifique, Sophia-Antipolis,
  France.}

\footnotetext[1]{School of Mathematical Sciences, Monash University, Melbourne, Australia
  (\texttt{david.wood@monash.edu}). Supported by the
Australian Research Council.}

\section{Introduction}

A major result in the theory of graph minors is that every graph $G$ with sufficiently large average degree $d(G)$ contains a complete graph $K_{t}$ as a minor. That is, a $K_{t}$ can be constructed from $G$ using \emph{vertex deletion}, \emph{edge deletion} and \emph{edge contraction}.
Let $$g(t) := \min\{D: \text{every graph G with } d(G) \geq D \text{ contains a } K_{t}\text{-minor}\}.$$ Mader \cite{fast1} showed that $g(t)$ is well-defined, and that $g(t) \leq 2^{t-2}$. Subsequently, Mader \cite{fast2} improved this bound to $g(t) \leq 16t\log_{2} t$, and later this was improved to $g(t) \in \Theta(t\sqrt{\log t})$ by Thomason \cite{fast3} and Kostochka \cite{fast4,fast5}, which is best possible. Thomason \cite{fast6} later determined the asymptotic constant for this bound.

This paper considers linear-time algorithms for finding a $K_{t}$-minor in a graph with high average degree. This question was first considered by Reed and Wood \cite{fast7} who gave a $O(n)$ time algorithm to find a $K_{t}$-minor in an $n$-vertex graph $G$ with $d(G) \geq 2^{t-2}$. We improve on this result by lowering the required bound on the average degree to within a constant factor of optimal:

\begin{theorem}
For all fixed $\epsilon > 0$ and fixed sufficiently large $t \geq t(\epsilon)$, there is a $O(n)$ time algorithm that, given an $n$-vertex graph $G$ with average degree $d(G) \geq (2+\epsilon)g(t)$, finds a $K_{t}$-minor in $G$. 
\end{theorem}


Reed and Wood used their algorithm mentioned above as a subroutine for finding separators in $H$-minor free graphs (also see \cite{fast11} for a related separator result). This result has subsequently been used in other algorithms for $H$-minor free graphs, in particular, shortest path algorithms by Tazari and M{\"u}ller-Hannemann \cite{fast8} and Wulff-Nilsen \cite{fast10}, and a maximum matching algorithm by Yuster and Zwick \cite{fast9}. The algorithm presented in this paper speeds up all these results (in terms of the dependence on $H$).

Finally, note that Robertson and Seymour \cite{fast12} describe a $O(n^{3})$ time algorithm that tests whether a given $n$-vertex graph contains a fixed graph $H$ as a minor. The time complexity was improved to $O(n^{2})$ by Kawarabayashi et al.~\cite{fast13}. Kawarabayashi and Reed have announced a $O(n \log n)$ time algorithm for this problem.

\renewcommand{\thefootnote}{\arabic{footnote}}

\section{Algorithm}
Given a vertex $v$ of a graph $G$, we denote by $\deg_{G}(v)$ and $N_{G}(v)$ the degree and neighbourhood of $v$ in $G$, respectively. We drop the subscript when $G$ is clear from the context.
Define a \emph{matching} $M \subseteq E(G)$ to be a set of edges such that no two edges in $M$ share an endpoint. Let $V(M)$ be the set of endpoints of the edges in $M$. An \emph{induced matching} in $G$ is a matching such that any two vertices $x,y$ of $V(M)$ are only adjacent in $G$ when $xy \in M$. Given a matching $M$ in $G$, let $G/M$ be the graph formed by contracting each edge of $M$ in $G$. 

We fix $\epsilon > 0$ and $t \geq 3$ such that $g(t) \geq \max\{t, \frac{2t}{\epsilon}\}$. We may assume $t \geq 3$ since finding a $K_{1}$- or $K_{2}$-minor is trivial, and that $g(t) \geq \max\{t, \frac{2t}{\epsilon}\}$ for sufficiently large $t$, since $g(t) \in \Theta(t\sqrt{\log t})$. Consider the following algorithm that takes as input a graph given as a list of vertices and a list of edges. The implicit output of the algorithm is the sequence of contractions and deletions that produce a $K_{t}$-minor.

\begin{algorithm}[H]
\caption{\textsc{FindMinor} (input: $n$-vertex graph $G$ with $d(G) \geq (2+\epsilon)g(t)$)}
\begin{algorithmic}[1]
\item[1:] Delete edges of $G$ so that $(2+\epsilon)g(t) \leq d(G) \leq (2+\epsilon)g(t)+1$.
\item[2:] Delete vertices of low degree so that the minimum degree $\delta(G) > \frac{1}{2}d(G)$.
\item[3:] Let $S := \{v \in V(G): \deg(v) \leq d(G)^{2}\}$, and let $B := \{v \in V(G): \deg(v) > d(G)^{2}\}$.
\newline [Note that $B$ is possibly empty, and that $S$ and $B$ partition $V(G)$.] 
\item[4:] Say an edge $vw \in E(G)$ is \emph{good} if $v,w \in S$ and $|N(v) \cap N(w)| \leq \frac{1}{2}(d(G)-2)$. Greedily construct a maximal matching $M$ of good edges.
\newline [Note that it is possible that no edges are good, in which case $M = \emptyset$.]
\item[5:] If $|M| > \frac{1}{8d(G)}n$, then greedily construct a maximal induced submatching $M'$ of $M$. That is, initialise $M':= \emptyset$ and $Q := M$, and repeat the following algorithm until $Q = \emptyset$: pick an edge $vw \in Q$, add $vw$ to $M'$, and delete from $Q$ the edge $vw$ and every edge with an endpoint adjacent to $v$ or $w$. 
\newline Let $G':=G/M'$. Run $\textsc{FindMinor}(G')$ and stop.
\item[6:] Now assume $|M| \leq \frac{1}{8d(G)}n$. Let $B' := B \cup V(M)$ and $S' := S-V(M)$.
\newline [Note that, similarly to Step 3, $S'$ and $B'$ partition $V(G)$.]
\item[7:] Greedily compute a maximal subset $A$ of $S'$ such that each vertex $u \in A$ is \emph{assigned} to a pair of vertices in $N(u) \cap B'$, and each pair of vertices in $B'$ has at most one vertex in $A$ assigned to it.
\item[8:] If $2|A| \geq d(G)|B'|$ and $B' \neq \emptyset$, then let $G'$ be the graph obtained from $G$ as follows: For each pair of distinct vertices $x,y \in B'$ with an assigned vertex $z \in A$, contract the edge $xz$. 
\newline Run $\textsc{FindMinor}(G'[B'])$ and stop.
\item[9:] Now assume $2|A| < d(G)|B'|$ or $B' = \emptyset$. Choose $v \in S'-A$. 
\newline [We prove below that $S'-A \neq \emptyset$. Since $v$ is not assigned, for every pair $x,y$ of vertices in $N(v) \cap B$ some vertex $z \in A$ is assigned to $x,y$.] 
\item[10:] If $|N(v) \cap B'| \geq t$, then let $G'$ be the graph obtained from $G$ as follows: For each pair of distinct vertices $x,y \in N(v) \cap B'$, if $z$ is the vertex in $A$ assigned to $x$ and $y$, then contract $xz$ into $x$ (so that the new vertex is in $B'$). Then $G'[N(v) \cap B'] \supseteq K_{t}$. Stop.
\item[11:] Otherwise let $G' := G[\{v\} \cup (N_{G}(v) \cap S')]$ and run an exhaustive search to find a $K_{t}$-minor in $G'$.
\newline [Below we prove that $d(G) \geq g(t)$ and $|V(G')| \leq d(G)^{2} + 1$.]
\end{algorithmic}
\end{algorithm}

\section{Correctness of Algorithm}
First, we prove that $\textsc{FindMinor}(G)$ does output a $K_{t}$-minor. Define $m:=|E(G)|$.
We must ensure the following: that $\textsc{FindMinor}$ finds a $K_{t}$-minor in Steps 5 and 8; that $S'-A \neq \emptyset$ in Step 9; that the graph constructed in Step 10 contains a $K_{t}$ subgraph; and that our exhaustive search in Step 11 finds a $K_{t}$-minor of $G$.

Consider Step 5. Assume that $\textsc{FindMinor}$ finds a $K_{t}$-minor in any graph $G'$ with $|V(G')| < n$ where $d(G') \geq (2+\epsilon)g(t)$. Consider the induced matching $M'$. Contracting any single edge $vw$ of $M'$ does not lower the average degree, as we only lose $|N(v) \cap N(w)|+1 \leq \frac{1}{2}d(G)$ edges and one vertex. Since the matching is induced, contracting every edge in $M'$ does not lower the average degree. Since $|M| > \frac{1}{8d(G)}n$, $M$ is not empty and $M'$ is not empty. Thus $d(G') \geq d(G) \geq (2+\epsilon)g(t)$ and $|V(G')| < |V(G)| =n$. Thus, by induction, running the algorithm on $G'$ finds a $K_{t}$-minor, and as such we find one for $G$.

If we recurse at Step 8, then $2|A| \geq d(G)|B'|$ and $B'\neq \emptyset$. Now $|V(G'[B'])|=|B'|$ and $|E(G'[B'])| \geq |A|$, since every assigned vertex corresponds to an edge of $G'[B']$. Thus $$d(G'[B']) = \frac{2|E(G'[B'])|}{|V(G'[B'])|} \geq \frac{2|A|}{|B'|} \geq d(G).$$ Also, $|V(G'[B'])|=|B'|<n$, since otherwise $A=S'=\emptyset$, contradicting $2|A| \geq d(G)|B'| > 0$. Hence, by assumption, the algorithm will find a $K_{t}$-minor in $G'[B']$. Thus the algorithm finds a $K_{t}$-minor for $G$. 

Now we show that $|S'| > |A|$ in Step 9. We have $2|A| < d(G)|B'|$ or $B' = \emptyset$. First consider the case when $2|A| < d(G)|B'|$. Note that $2m=d(G)n$, and that $d(G)^{2}|B| < \sum\nolimits_{v \in B} \deg(v) \leq 2m$, and so $|B| < \frac{2m}{d(G)^{2}} = \frac{1}{d(G)}n$. 
Now $|S'|=|S|-2|M| \geq |S| - \frac{1}{4d(G)}n$ by Step 6. Thus, $$|S'| \geq |S| - \frac{1}{4d(G)}n = (n - |B|) - \frac{1}{4d(G)}n > n - \frac{1}{d(G)}n - \frac{1}{4d(G)}n = \frac{4d(G)-5}{4d(G)}n.$$
By Step 9 and Step 6, $$|A| < \frac{d(G)}{2}|B'| = \frac{d(G)}{2}(|B| + 2|M|) < \frac{d(G)}{2}\left(\frac{1}{d(G)}n + \frac{1}{4d(G)}n\right) = \frac{5}{8}n.$$
Thus, if $|S'| \leq |A|$ then $\frac{4d(G)-5}{4d(G)}n < \frac{5}{8}n$, so $3d(G) < 10$, which is a contradiction since $d(G) \geq (2+\epsilon)g(t) > 2g(3) = 4$. (We have $g(t) \geq g(3) = 2$, since $g(t)$ is non-decreasing.) Hence, $|S'| > |A|$. Now consider the case that $B'= \emptyset$. Then $|S'|=n$ and $A=\emptyset$, since the vertices of $A$ are assigned to pairs of vertices in $B'$. Hence $|S'| > |A|$.

Now consider Step 10. $G'[N(v) \cap B']$ has at least $t$ vertices by assumption. Each pair of distinct vertices $x,y$ in $N(v) \cap B'$ has an assigned vertex in $A$, as otherwise $v$ would have been assigned to $x$ and $y$. Hence the vertex $z$ exists, and $x$ and $y$ are adjacent after contracting $xz$. Therefore all pairs of vertices in $N(v) \cap B'$ become adjacent, and $G'[N(v) \cap B']$ is a complete graph, and we have found our $K_{t}$-minor in $G$. 

Finally consider Step 11. $G'$ is an induced subgraph of $G$, and so if we can find $K_{t}$ as a minor in $G'$, we have a $K_{t}$-minor in $G$. We use an exhaustive search, so all we need to ensure is that $G'$ does have a $K_{t}$-minor. Thus, we simply need to ensure that $d(G') \geq g(t)$. By Step 1 and Step 2, $\deg_{G}(v) > \frac{1}{2}d(G) \geq \frac{\epsilon}{2}g(t) \geq t$. Since Step 10 was not applicable, $v$ has at most $t-1$ neighbours in $B'$. Thus $v$ has some neighbour in $S'$. Let $w$ be a vertex of $G'-v$. Thus $vw$ is an edge and $v,w \in S'$. Since neither $v$ nor $w$ was matched by $M$, and since $M$ is maximal, $vw$ is not good. Since $v,w \in S' \subseteq S$, this means that $|N(v) \cap N(w)| > \frac{1}{2}(d(G)-2)$. As $v$ has at most $t-1$ neighbours in $B'$, we have $|N(v) \cap N(w) \cap S'| > \frac{1}{2}(d(G)-2) - (t-1)$. Every common neighbour of $v$ and $w$ in $S'$ is a neighbour of $w$ in $G'$, by definition, so $\deg_{G'}(w) > \frac{1}{2}(d(G)-2) - (t-1)$. Since $v$ is dominant in $G'$, we have $d(G') \geq \frac{1}{2}(d(G)-2) - (t-1)$, which is at least $g(t)$ as required since $d(G) \geq (2+\epsilon)g(t)$ and $\epsilon g(t) \geq 2t$. 

\section{Time Complexity}
Now that we have shown that $\textsc{FindMinor}$ will output a $K_{t}$-minor, we must ensure it does so in $O(n)$ time (for fixed $t$ and $\epsilon$). 

First, suppose $\textsc{FindMinor}$ runs without recursing.
Recall that our input graph $G$ is given as a list of vertices and a list of edges, from which we will construct adjacency lists as it is read in. Since our goal in Step 1 is to ensure that $m \leq \frac{1}{2}((2+\epsilon)g(t)+1)n$, we can do this by taking, at most, the first $\frac{1}{2}((2+\epsilon)g(t)+1)n$ edges, and ignoring the rest. This can be done in $O(n)$ time, and from now on we may assume that $m \in O(n)$.
In Step 2, since we are only deleting vertices of bounded degree, this can be done in $O(n)$ time.
Clearly, Steps 3, 6 and 9 can be implemented in $O(n)$ time. By definition, the degree of any vertex in $S$ or $S'$ is at most $((2+\epsilon)g(t)+1)^{2}$. Hence Steps 4, 5, 7, 8 and 10 take $O(n)$ time. Finally, for Step 11  note that $|V(G')| \leq d(G)^{2} + 1$, so exhaustive search runs in $O(1)$ time for fixed $t$.
Hence the algorithm without recursion runs in $O(n)$ time.

Should $\textsc{FindMinor}$ recurse, we need to ensure that the order of the graph we recurse on is a constant factor less than $n$. Then the overall time complexity is $O(n)$ (by considering the sum of a geometric series). 
In Step 5, the endpoints of edges in $M$ have degree at most $d(G)^{2}$, and thus $|M'| \geq \frac{1}{2d(G)^2}|M| \geq \frac{1}{16d(G)^3}n$. This ensures that $|V(G')| \leq (1-\frac{1}{16d(G)^{3}})n$, as desired. In Step 8, the order of $G'[B']$ is at most $\frac{2|A|}{d(G)} \leq \frac{2n}{d(G)}$. Hence it follows that the overall time complexity is $O(n)$.

\bibliographystyle{plain}
\bibliography{fastbib}

\end{document}